\documentclass{article}

\usepackage{PRIMEarxiv}
\usepackage{amsmath}
\usepackage{subcaption} 

\usepackage[linesnumbered,ruled,vlined]{algorithm2e}
\usepackage{float}
\usepackage[utf8]{inputenc} % allow utf-8 input
\usepackage[T1]{fontenc}    % use 8-bit T1 fonts
\usepackage{hyperref}       % hyperlinks
\usepackage{url}            % simple URL typesetting
\usepackage{booktabs}       % professional-quality tables
\usepackage{amsfonts}       % blackboard math symbols
\usepackage{nicefrac}       % compact symbols for 1/2, etc.
\usepackage{microtype}      % microtypography
\usepackage{lipsum}
\usepackage{fancyhdr}       % header
\usepackage{graphicx}       % graphics
\graphicspath{{media/}}     % organize your images and other figures under media/ folder

\title{Newton's Method Applied to Nonlinear Boundary Value Problems: A Numerical Approach
}

\author{
  Rusdrael Antony de Araújo Freire \\
  Federal University of Rio Grande do Norte \\
  Caicó- RN\\
  \texttt{antony.freire.086@edu.ufrn.br} \\
   \And
  Francisco Márcio Barboza \\
  Department of Computing and Technology \\
  Federal University of Rio Grande do Norte \\
  Caicó- RN\\
  \texttt{marcio.barboza@ufrn.br} \\
   \AND
   Márcio Matheus de Lima Barboza \\
 Federal University of Rio Grande do Norte \\
  Natal- RN\\
  \texttt{marcio.barboza.024@edu.ufrn.br} \\
}

\begin{document}
\maketitle

\begin{abstract}
This work investigates the application of the Newton's method for the numerical solution of a nonlinear boundary value problem formulated through an ordinary differential equation (ODE). Nonlinear ODEs arise in various mathematical modeling contexts, where an exact solution is often unfeasible due to the intrinsic complexity of these equations. Thus, a numerical approach is employed, using Newton's method to solve the system resulting from the discretization of the original problem. The procedure involves the iterative formulation of the method, which enables the approximation of solutions and the evaluation of convergence with respect to the problem parameters. The results demonstrate that Newton's method provides a robust and efficient solution, highlighting its applicability to complex boundary value problems and reinforcing its relevance for the numerical analysis of nonlinear systems. It is concluded that the methodology discussed is suitable for solving a wide range of boundary value problems, ensuring precision and stability in the results.
\end{abstract}

% keywords can be removed

\section{Introduction}
Ordinary differential equations (ODEs) play a central role in mathematical modeling, being used in fields ranging from physics and engineering to biology and economics \cite{Boyce2017,Smith2009}. In particular, boundary value problems for ODEs appear in various applications, such as the description of heat transfer phenomena, fluid mechanics, and population growth \cite{Tenenbaum1963}. However, the exact solution of these equations is not always possible, especially when the equation is nonlinear, which requires numerical methods to obtain approximate solutions \cite{remani2013numerical}.

Among numerical methods, Newton's method has been widely used due to its efficiency and rapid convergence, especially in nonlinear systems \cite{Burden2011}. In the context of nonlinear ODEs, it becomes a powerful tool for solving boundary value problems by transforming the initial problem into a system of nonlinear equations, whose solution is iteratively approximated \cite{Cheney2012}. This approach is particularly relevant in situations where the behavior of the solution is complex and requires refinement of successive approximations.

The present work aims to apply Newton's method to solve a nonlinear boundary value problem, based on the formulation proposed by \cite{remani2013numerical}. The differential equation addressed represents a typical example of a nonlinear problem, and the boundary conditions impose additional challenges that require robust numerical strategies. Using an interval discretization and iterative formulations, this study presents an approximate solution to the proposed problem, along with a detailed analysis of the convergence and accuracy of the results.

Thus, the objective of this work is to demonstrate the effectiveness of Newton's method in solving nonlinear boundary value problems, providing a practical approach and a theoretical basis that can be applied to a wide range of similar problems. The analysis of the results obtained illustrates the potential of this method in providing stable and reliable numerical solutions for complex ODEs, reinforcing its applicability and importance in the field of numerical analysis.

\section{Newton's Method}

Newton's method is a highly efficient iterative technique for solving nonlinear equations and finding the roots of such functions. It is applicable to both single-variable functions and systems with multiple variables. The main idea is to replace the nonlinear function with a local linear approximation to obtain a better estimate of the root.

\subsection{Newton's Method in One Dimension}

For a single-variable function \( f(x) \), Newton's method uses the tangent line to the curve of the function \( f \) at an initial point \( x_0 \), which should be relatively close to the true root \( \alpha \). The iterative formula to obtain the next approximation \( x_{n+1} \) is given by:
\[
x_{n+1} = x_n - \frac{f(x_n)}{f'(x_n)}.
\]
In this case, \( f'(x_n) \) represents the derivative of \( f(x) \) evaluated at \( x_n \). This process approximates each subsequent point through the tangent, allowing for rapid convergence towards the root \( \alpha \), provided that \( x_0 \) is sufficiently close to it and that \( f'(x_n) \neq 0 \) in all iterations.

\begin{figure}[H]
    \centering
    \includegraphics[width =15cm]{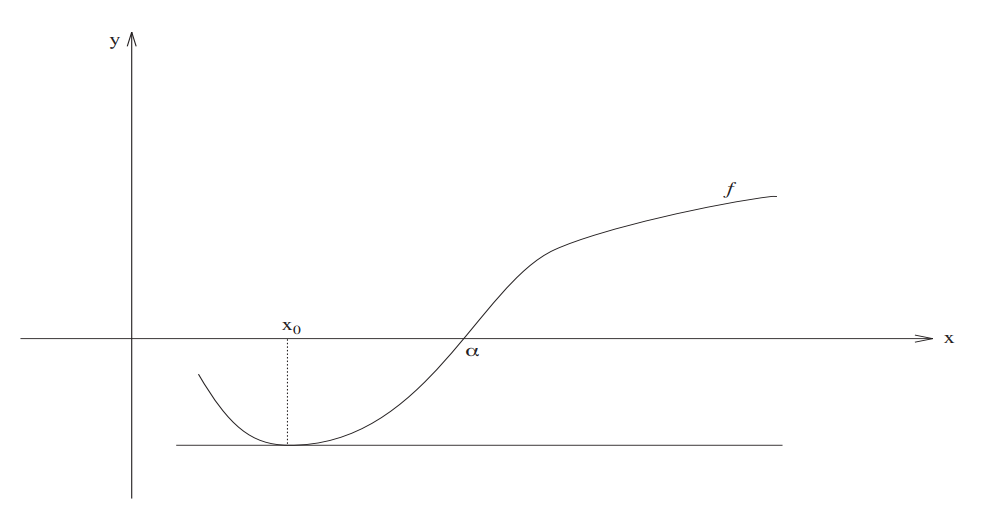}
    \caption{Failure of the Newton Method due to a poor initial point \cite{kharab2011introduction}.}
    \label{fig:newton-grafico-falha}
\end{figure}

If \( x_0 \) is far from \( \alpha \) or if \( f'(x_0) \) is close to zero, the method may diverge, as illustrated in Figure \ref{fig:newton-grafico-falha}. The method has a quadratic convergence rate when \( x_0 \) is chosen appropriately \cite{kharab2011introduction}.

\subsection{Newton's Method for Nonlinear Systems}

The concept of Newton's method can be extended to solve systems of nonlinear equations, where we seek simultaneous solutions for multiple functions of multiple variables. Consider a system with two nonlinear equations:
\[
f_1(x_1, x_2) = 0,
\]
\[
f_2(x_1, x_2) = 0.
\]
Newton's method for this system begins with an initial estimate \((x_1^{(0)}, x_2^{(0)})\). In each iteration, each function is approximated by its tangent plane, and the solution of the resulting linear system provides the next approximation of the root. The vector formulation to obtain this new estimate is:
\[
J(\mathbf{X}^{(i)}) \Delta \mathbf{X}^{(i)} = -\mathbf{F}(\mathbf{X}^{(i)}),
\]
where \( J(\mathbf{X}^{(i)}) \) is the Jacobian matrix of the system, composed of the partial derivatives of the functions with respect to the variables. In the case of two variables, the Jacobian matrix is:
\[
J(\mathbf{X}^{(i)}) = \begin{bmatrix} \frac{\partial f_1}{\partial x_1} & \frac{\partial f_1}{\partial x_2} \\ \frac{\partial f_2}{\partial x_1} & \frac{\partial f_2}{\partial x_2} \end{bmatrix}.
\]
After solving this system for \( \Delta \mathbf{X}^{(i)} = (\Delta x_1^{(i)}, \Delta x_2^{(i)}) \), we update the solution:
\[
x_1^{(i+1)} = x_1^{(i)} + \Delta x_1^{(i)}, \quad x_2^{(i+1)} = x_2^{(i)} + \Delta x_2^{(i)}.
\]

This process is repeated until the vector \( \mathbf{F}(\mathbf{X}^{(i)}) \) is close to zero or the variations in the variables are sufficiently small, indicating convergence. For a system with \( n \) equations and \( n \) variables, the generalized notation is:
\[
J(\mathbf{X}) \Delta \mathbf{X} = -\mathbf{F}(\mathbf{X}),
\]
where \( J(\mathbf{X}) \) is the \( n \times n \) Jacobian matrix. The solution update is:
\[
\mathbf{X}^{(i+1)} = \mathbf{X}^{(i)} + \Delta \mathbf{X}^{(i)}.
\]

This extension of Newton's method to nonlinear systems is powerful, but it requires that the Jacobian matrix be nonsingular and that the initial solution choice be close to the true value to ensure quadratic convergence \cite{kharab2011introduction}.

As discussed, Newton's method involves iteratively solving a system of nonlinear equations using linear approximations. The process can be described by the following pseudocode:

\begin{algorithm}[H]
\caption{Pseudocode for Newton's Method for Nonlinear Systems}\label{alg:newton}
\KwIn{$F$: The nonlinear function we want to solve}
\KwIn{$JF$: The function that calculates the Jacobian of $F$ (can be empty)}
\KwIn{$\mathbf{X}_0$: The initial guess vector for the solution}
\KwIn{$tol$: The tolerance for the stopping criterion}
\KwIn{$maxit$: The maximum number of allowed iterations}
\SetKwFunction{FMain}{Newton's Method}
\SetKwProg{Fn}{Function}{:}{}
\Fn{\FMain}{
    $\mathbf{X} \gets \mathbf{X}_0$ \;
    $iter \gets 1$ \;
    $h \gets 10^{-6}$ \tcp{Small step for finite differences}
    \While{$iter \leq maxit$}{
        \If{$JF$ is empty}{
            $J \gets$ calculate Jacobian numerically using $F$, $\mathbf{X}$, and $h$ \;
        }
        \Else{
            $J \gets$ call function $JF$ with $\mathbf{X}$ \;
        }
        $\Delta \mathbf{X} \gets - J^{-1} F(\mathbf{X})$ \tcp{Newton step}
        $\mathbf{X}_n \gets \mathbf{X} + \Delta \mathbf{X}$ \tcp{Update the solution}
        $err \gets$ maximum absolute value of $(\mathbf{X}_n - \mathbf{X})$ \tcp{Calculate the error}
        \If{$err \leq tol$}{
            \KwRet{$\mathbf{X}_n$} \tcp{Convergence achieved}
        }
        \Else{
            $\mathbf{X} \gets \mathbf{X}_n$ \;
        }
        $iter \gets iter + 1$ \;
    }
    \textbf{Error:} Newton's method did not converge \;
    \KwRet{Last calculated solution $\mathbf{X}_n$}
}
\end{algorithm}

\section{Problem Definition}
This study is based on the work of \cite{remani2013numerical}, which explores the solution of boundary value problems for nonlinear ordinary differential equations. These problems are relevant in various applications, from modeling physical phenomena to engineering and material science. The complexity of nonlinear problems, particularly in the context of fixed boundary conditions, requires robust numerical approaches to find reliable approximate solutions.

In this paper, we consider the following boundary value problem involving a nonlinear ordinary differential equation:
\begin{equation}
y'' = \frac{1}{8}(32 + 2x^3 - yy'), \quad 1 \leq x \leq 3
\end{equation}
The imposed boundary conditions are:
\begin{equation}
y(1) = 17, \quad y(3) = 14.333333
\end{equation}
This problem describes the relationship between the dependent variable \( y \), its first derivative \( y' \), and its second derivative \( y'' \), subject to specific values of \( y \) at \( x = 1 \) and \( x = 3 \). The nonlinear function of the derivative and the function itself makes the analytical solution unattainable, motivating the application of numerical methods for a satisfactory approximation.

To solve the problem, we choose to discretize the interval \([1, 3]\) using a step size \( h = 0.1 \), which generates a sequence of uniformly spaced points \( x_i \). This approach allows for the formulation of a system of nonlinear equations from the discretization of the differential equation, which will then be solved iteratively using Newton's method. The choice of \( h = 0.1 \) is a compromise between accuracy and computational cost, ensuring that the resulting system is sufficiently detailed to capture the solution's dynamics while remaining computationally feasible.

\section{Numerical Method}

The numerical solution of the proposed boundary value problem is carried out using the interval discretization method and an iterative approximation technique. The steps of the adopted method are detailed below.

\subsection{Step 1: Interval Discretization}

The first step consists of discretizing the interval \([1, 3]\), where the value of \( x \) varies between 1 and 3. To do this, we divide the interval into \( N + 1 = 20 \) subintervals of equal width, with a step size \( h = 0.1 \). This means that the distance between two consecutive points \( x_i \) and \( x_{i+1} \) will be 0.1. Each value \( x_i \) is calculated by:
\begin{equation}
x_i = a + i \times h, \quad i = 0, 1, 2, \ldots, 20\end{equation}
where \( a = 1 \) is the lower bound of the interval, and \( h = 0.1 \) is the step size. The resulting values of \( x_i \) for each \( i \) are presented in Table \ref{table:xi_values}, which illustrates the discretization of the interval.

\begin{table}[H]
\centering
\begin{tabular}{|c|c|}
\hline
$i$ & $x_i$ \\
\hline
0 & 1.0 \\
1 & 1.1 \\
2 & 1.2 \\
3 & 1.3 \\
4 & 1.4 \\
5 & 1.5 \\
6 & 1.6 \\
7 & 1.7 \\
8 & 1.8 \\
9 & 1.9 \\
10 & 2.0 \\
11 & 2.1 \\
12 & 2.2 \\
13 & 2.3 \\
14 & 2.4 \\
15 & 2.5 \\
16 & 2.6 \\
17 & 2.7 \\
18 & 2.8 \\
19 & 2.9 \\
20 & 3.0 \\
\hline
\end{tabular}
\caption{Values of $x_i$ for the interval $[1, 3]$ with $h = 0.1$}
\label{table:xi_values}
\end{table}

\subsection{Step 2: Definition of Boundary Conditions}

The boundary conditions for the problem are provided as values of \( y \) at the interval's endpoints. In this case, we have:
\[
w_0 = 17 \quad \text{and} \quad w_{20} = 14.333333
\]
These values are used as starting points for the numerical solution and will be applied in the iteration process.

\subsection{Step 3: Initial Approximation}

To start the iterative process, an initial approximation \( w^{(0)} \) for the values of \( y \) at all discretized points \( x_i \) must be defined. The initial approximation is given by:

\begin{equation}
w^{(0)} = \begin{pmatrix}
15.6666 \\
15.6666 \\
15.6666 \\ 
15.6666 \\ 
15.6666 \\
15.6666 \\
15.6666 \\
15.6666 \\
15.6666 \\
15.6666 \\
15.6666 \\
15.6666 \\
15.6666 \\
15.6666 \\
15.6666 \\
15.6666 \\
15.6666 \\
15.6666 \\
15.6666
\end{pmatrix}^t
\end{equation}

The initial approximation for the values \( w_i \), with \( i = 1, 2, \ldots, 19 \), was obtained from the average of the boundary values \( w_0 = 17 \) and \( w_{20} = 14.33 \), as defined by the boundary conditions.

\subsection{Step 4: Definition of the Nonlinear System}

The nonlinear system that describes the problem consists of 19 equations, which relate the values of \( w_i \) to their respective approximations. Each equation in the system \( F(w) = 0 \) has the general form:

\begin{equation}
F(w) = \begin{pmatrix}
2w_1 - w_2 + 0.01 \left( 4 + 0.33275 + \frac{w_1(w_2 - 17)}{1.6} \right) - 17 &= 0 \\
-w_1 + 2w_2 - w_3 + 0.01 \left( 4 + 0.432 + \frac{w_2(w_3 - w_1)}{1.6} \right) &= 0 \\
-w_2 + 2w_3 - w_4 + 0.01 \left( 4 + 0.5495 + \frac{w_3(w_4 - w_2)}{1.6} \right) &= 0 \\
-w_3 + 2w_4 - w_5 + 0.01 \left( 4 + 0.686 + \frac{w_4(w_5 - w_3)}{1.6} \right) &= 0 \\
-w_4 + 2w_5 - w_6 + 0.01 \left( 4 + 0.84375 + \frac{w_5(w_6 - w_4)}{1.6} \right) &= 0 \\
-w_5 + 2w_6 - w_7 + 0.01 \left( 4 + 1.024 + \frac{w_6(w_7 - w_5)}{1.6} \right) &= 0 \\
-w_6 + 2w_7 - w_8 + 0.01 \left( 4 + 1.22825 + \frac{w_7(w_8 - w_6)}{1.6} \right) &= 0 \\
-w_7 + 2w_8 - w_9 + 0.01 \left( 4 + 1.458 + \frac{w_8(w_9 - w_7)}{1.6} \right) &= 0 \\
-w_8 + 2w_9 - w_{10} + 0.01 \left( 4 + 1.71475 + \frac{w_9(w_{10} - w_8)}{1.6} \right) &= 0 \\
-w_9 + 2w_{10} - w_{11} + 0.01 \left( 4 + 2 + \frac{w_{10}(w_{11} - w_9)}{1.6} \right) &= 0 \\
-w_{10} + 2w_{11} - w_{12} + 0.01 \left( 4 + 2.31525 + \frac{w_{11}(w_{12} - w_{10})}{1.6} \right) &= 0 \\
-w_{11} + 2w_{12} - w_{13} + 0.01 \left( 4 + 2.662 + \frac{w_{12}(w_{13} - w_{11})}{1.6} \right) &= 0 \\
-w_{12} + 2w_{13} - w_{14} + 0.01 \left( 4 + 3.04175 + \frac{w_{13}(w_{14} - w_{12})}{1.6} \right) &= 0 \\
-w_{13} + 2w_{14} - w_{15} + 0.01 \left( 4 + 3.456 + \frac{w_{14}(w_{15} - w_{13})}{1.6} \right) &= 0 \\
-w_{14} + 2w_{15} - w_{16} + 0.01 \left( 4 + 3.90625 + \frac{w_{15}(w_{16} - w_{14})}{1.6} \right) &= 0 \\
-w_{15} + 2w_{16} - w_{17} + 0.01 \left( 4 + 4.394 + \frac{w_{16}(w_{17} - w_{15})}{1.6} \right) &= 0 \\
-w_{16} + 2w_{17} - w_{18} + 0.01 \left( 4 + 4.92075 + \frac{w_{17}(w_{18} - w_{16})}{1.6} \right) &= 0 \\
-w_{17} + 2w_{18} - w_{19} + 0.01 \left( 4 + 5.488 + \frac{w_{18}(w_{19} - w_{17})}{1.6} \right) &= 0 \\
-w_{18} + 2w_{19} + 0.01 \left( 4 + 6.09725 + \frac{w_{19}(14.333333 - w_{18})}{1.6} \right) - 14.333333 &= 0 \\
\end{pmatrix}
\end{equation}

This system is derived from the discretization of the original differential equation, considering interactions between adjacent points \( w_i \) and the boundary conditions. The objective is to find a numerical solution that satisfies this system of equations.

\subsection{Step 5: Calculation of the Jacobian \( J(w) \)}

The Jacobian \( J(w) \) of the nonlinear system is a matrix that describes the partial derivatives of \( F(w) \) with respect to each component \( w_i \). This matrix is used in the Newton method to update the solution approximations. The Jacobian matrix \( J(w) \) is structured as follows:

\[
J(w) =
\begin{bmatrix}
2 + 0.01 (\frac{w_2 - 17}{1.6}) & -1 + 0.05 (\frac{1}{8}w_1) & 0 & \cdots & \cdots & \cdots & \cdots & 0 \\
-1 - 0.05 (\frac{1}{8}w_2) & 2 + 0.01 (\frac{w_3 - w_1}{1.6}) & -1 + 0.05 (\frac{1}{8}w_2) & 0 & \cdots & \cdots & \cdots & 0 \\
0 & a_{32} & a_{33} & a_{34} & 0 & \cdots & \cdots & 0 \\
\vdots & 0 & \ddots & \ddots & \ddots & 0 & \cdots & 0 \\
\vdots & \vdots & 0 & \ddots & \ddots & \ddots & 0 & 0 \\
\vdots & \vdots & \vdots & 0 & \ddots & \ddots & \ddots & 0 \\
\vdots & \vdots & \vdots & \vdots & 0 & \ddots & \ddots & a_{i-1j} \\
0 & 0 & 0 & 0 & 0 & 0 & a_{i,j-1} & a_{ij} \\
\end{bmatrix}
\]

The Jacobian matrix is essential to the iterative process, as it is used to calculate the corrections needed for the \( w_i \) values in each iteration.

\section{Results}

The iterative method was applied successively until the difference between consecutive iterations became sufficiently small, indicating the solution's convergence. The convergence criterion was set as \( ||w^{(k)} - w^{(k-1)}|| \leq \epsilon \), where \( \epsilon \) is a predefined tolerance value, signifying that the solution achieved satisfactory accuracy. Table \ref{table:results} displays the \( w^{(k)} \) values for each iteration \( k \), until convergence was reached, with the condition \( ||w^{(4)} - w^{(3)}|| = 0 \) demonstrating that the \( w \) values in successive iterations stabilized.

\begin{table}[H]
\centering
\begin{tabular}{|c|c|c|c|c|c|c|}
\hline
$x_i$ & $w_i$ & $\mathbf{w}^{(0)}$ & $\mathbf{w}^{(1)}$ & $\mathbf{w}^{(2)}$ & $\mathbf{w}^{(3)}$ & $\mathbf{w}^{(4)}$ \\
\hline
1.0 & $w_0$    & 17.0000 & 17.0000 & 17.0000 & 17.0000 & 17.0000 \\
1.1 & $w_1$    & 15.6666 & 16.7657 & 16.7605 & 16.7605 & 16.7605 \\
1.2 & $w_2$    & 15.6666 & 16.5183 & 16.5135 & 16.5134 & 16.5134 \\
1.3 & $w_3$    & 15.6666 & 16.2664 & 16.2589 & 16.2589 & 16.2589 \\
1.4 & $w_4$    & 15.6666 & 16.0102 & 15.9974 & 15.9974 & 15.9974 \\
1.5 & $w_5$    & 15.6666 & 15.7504 & 15.7299 & 15.7298 & 15.7298 \\
1.6 & $w_6$    & 15.6666 & 15.4879 & 15.4578 & 15.4577 & 15.4577 \\
1.7 & $w_7$    & 15.6666 & 15.2240 & 15.1830 & 15.1829 & 15.1829 \\
1.8 & $w_8$    & 15.6666 & 14.9609 & 14.9084 & 14.9083 & 14.9083 \\
1.9 & $w_9$    & 15.6666 & 14.7011 & 14.6376 & 14.6375 & 14.6375 \\
2.0 & $w_{10}$ & 15.6666 & 14.4483 & 14.3751 & 14.3750 & 14.3750 \\
2.1 & $w_{11}$ & 15.6666 & 14.2070 & 14.1267 & 14.1266 & 14.1266 \\
2.2 & $w_{12}$ & 15.6666 & 13.9835 & 13.8994 & 13.8993 & 13.8993 \\
2.3 & $w_{13}$ & 15.6666 & 13.7852 & 13.7019 & 13.7018 & 13.7018 \\
2.4 & $w_{14}$ & 15.6666 & 13.6220 & 13.5444 & 13.5443 & 13.5443 \\
2.5 & $w_{15}$ & 15.6666 & 13.5060 & 13.4392 & 13.4391 & 13.4391 \\
2.6 & $w_{16}$ & 15.6666 & 13.4525 & 13.4010 & 13.4010 & 13.4010 \\
2.7 & $w_{17}$ & 15.6666 & 13.4804 & 13.4475 & 13.4475 & 13.4475 \\
2.8 & $w_{18}$ & 15.6666 & 13.6132 & 13.5999 & 13.5999 & 13.5999 \\
2.9 & $w_{19}$ & 15.6666 & 13.8801 & 13.8843 & 13.8843 & 13.8843 \\
3.0 & $w_{20}$ & 14.3333 & 14.3333 & 14.3333 & 14.3333 & 14.3333 \\
\hline
\end{tabular}
\caption{Iteration results for the value of $w_i$}
\label{table:results}
\end{table}

\section{Conclusion}

After applying the iterative method with the Jacobian matrix, an approximate solution to the boundary value problem was obtained through successive iterations. The approximate values for \( w_i \) across the discretized interval are given by the sequence:

\begin{equation}
\overline{\mathbf{w}} = \begin{pmatrix}
17.0000, \\
16.7605, \\
16.5134, \\
16.2589, \\
15.9974, \\
15.7298, \\
15.4577, \\
15.1829, \\
14.9083, \\
14.6375, \\
14.3750, \\
14.1266, \\
13.8993, \\
13.7018, \\
13.5443, \\
13.4391, \\
13.4010, \\
13.4475, \\
13.5999, \\
13.8843, \\
14.3333
\end{pmatrix}^t
\end{equation}

Each value \( w_i \) represents an approximation of \( y(x_i) \), with \( x_i = a + i \times h \) and \( i = 0, 1, 2, \ldots, N+1 \), where \( h \) is the spacing between discretized points. The Newton method, when applied to the system of discretized differential equations, effectively converged to a solution for this boundary value problem, as demonstrated by the stabilization of \( w_i \) values in the final iterations.

Below, we present a graph illustrating the evolution of the vector \( \overline{\mathbf{w}} \) over the iterations. This graph visually demonstrates how the values \( w_i \) stabilize, reflecting the convergence of the method to the approximate solution.

\begin{figure}[H]
    \centering
    \includegraphics[width=0.5\linewidth]{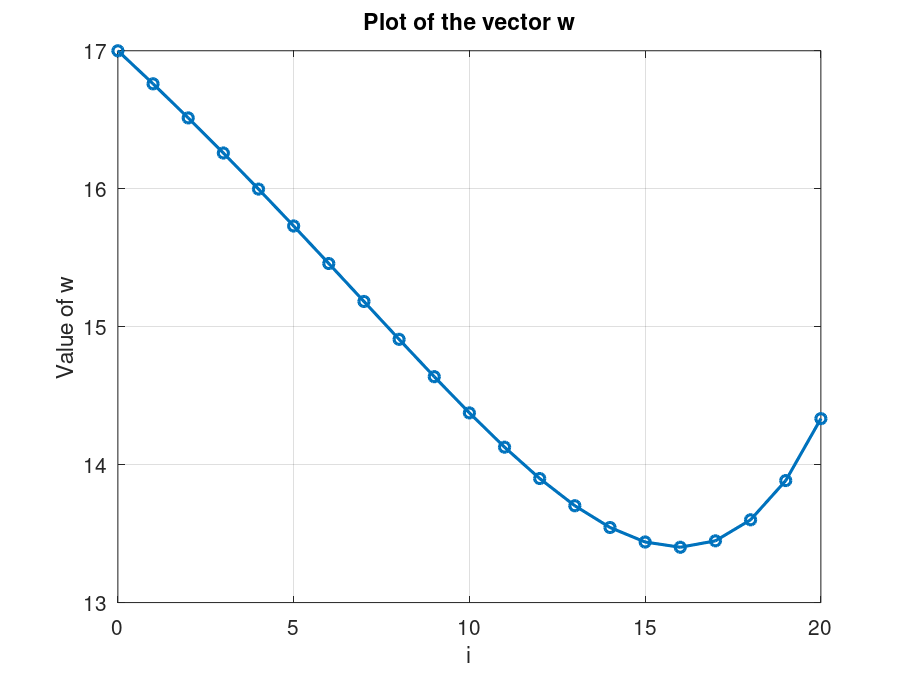}
        \caption{Plot of the vector \( \overline{\mathbf{w}} \) over the discretized interval.}
    \label{fig:problem_plot}
\end{figure}

%Bibliography
\bibliographystyle{unsrt}  
\bibliography{references}

\end{document}